\documentclass[11pt,reqno]{amsart}
\usepackage{enumerate}

\usepackage{color}

\newtheorem{theorem}{Theorem}
\newtheorem{proposition}{Proposition}
\newtheorem{corollary}[proposition]{Corollary}
\newtheorem{remark}[proposition]{Remark}
\newtheorem{lemma}[proposition]{Lemma}
\newfont{\bb}{msbm10 at 12pt}

\def\pf{{\textit {Proof :} }}

\def\C{{\mathbb C}}

\def\S{\hbox{\bb S}}

\def\bS{{\mathbb S}}
\def\what{\widehat}
\def\E{\mathcal{E}}
\def\SM{{\mathbb{S} M}}
\def\<{\langle}     
\def\>{\rangle}     

\def\mby{m_{BY}(\Sigma)}
\def\ml{m_{L}(\Sigma)}
\def\mhm{m(\Sigma)}
\def\mg{\mathcal{M}(\Sigma)}

\newcommand{\bal}{\begin{align}}      \newcommand{\eal}{\end{align}}
\newcommand{\ba}{\begin{array}}      \newcommand{\ea}{\end{array}}
\newcommand{\bc}{\begin{center}}     \newcommand{\ec}{\end{center}}
\newcommand{\be}{\begin{enumerate}}  \newcommand{\ee}{\end{enumerate}}
\newcommand{\beQ}{\begin{eqnarray*}} \newcommand{\eeQ}{\end{eqnarray*}}
\newcommand{\bi}{\begin{itemize}}    \newcommand{\ei}{\end{itemize}}
\newcommand{\bt}{\begin{tabular}}    \newcommand{\et}{\end{tabular}}
\newcommand{\bdm}{\begin{displaymath}} \newcommand{\edm}{\end{displaymath}}

    \newcommand{\sm}{\bS\!\!\!/\,\!}
\newcommand{\D}{D\!\!\!\!/\,}
\newcommand{\nb}{\nabla\!\!\!\!/\,}
\newcommand{\mult}{\gamma\!\!\!/}

\def\qed{\hfill{q.e.d.}\smallskip\smallskip}

\begin{document}

\title[On a Liu-Yau type Inequality]{On a Liu-Yau type Inequality for Surfaces} 

\author{Oussama Hijazi}
\address[Oussama Hijazi]{Institut {\'E}lie Cartan de Lorraine,
Universit{\'e} de Lorraine, Nancy,
B.P. 239,
54506 Vand\oe uvre-L{\`e}s-Nancy Cedex, France.}
\email{Oussama.Hijazi@univ-lorraine.fr}

\author{Sebasti{\'a}n Montiel}
\address[Sebasti{\'a}n  Montiel]{Departamento de Geometr{\'\i}a y Topolog{\'\i}a,
Universidad de Granada,
18071 Granada,  Spain.}
\email{smontiel@ugr.es}

\author{Simon Raulot}
\address[Simon Raulot]{Laboratoire de Math\'ematiques R. Salem
UMR $6085$ CNRS-Universit\'e de Rouen
Avenue de l'Universit\'e, BP.$12$
Technop\^ole du Madrillet
$76801$ Saint-\'Etienne-du-Rouvray, France.}
\email{simon.raulot@univ-rouen.fr}

\begin{abstract}
Let $\Omega$ be a compact and mean-convex domain with smooth boundary $\Sigma:=\partial\Omega$, in an initial data set $(M^3,g,K)$, which has no apparent horizon in its interior. If $\Sigma$ is spacelike in a spacetime $(\E^4,g_\E)$ with spacelike mean curvature vector $\mathcal{H}$ such that $\Sigma$ admits an isometric and isospin immersion into $\mathbb{R}^3$ with mean curvature $H_0$, then: 
\begin{eqnarray*}
\int_{\Sigma}|\mathcal{H}|d\Sigma\leq\int_{\Sigma}\frac{H_0^2}{|\mathcal{H}|}d\Sigma.
\end{eqnarray*}
If equality occurs, we prove that there exists a local isometric immersion of $\Omega$ in $\mathbb{R}^{3,1}$ (the Minkowski spacetime) with second fundamental form given by $K$. In Theorem \ref{liu-yau-minkowski}, we also examine, under weaker conditions, the case where the spacetime is the $(n+2)$-dimensional Minkowski space $\mathbb{R}^{n+1,1}$ and establish a stronger rigidity result.
\end{abstract}



\thanks{The second author was partially  
supported by a Spanish MEC-FEDER grant No. MTM2007-61775}

\date{\today}   

\maketitle 
\pagenumbering{arabic}
 

\section{Introduction}


Let $(\E^4,g_{\E})$ be a spacetime satisfying the Einstein field equations, that is $(\E^4,g_{\E})$ is a $4$-dimensional time-oriented Lorentzian manifold such that
\begin{eqnarray*}
Ric_\E-\frac{1}{2} R_\E g_{\E}=\mathcal{T},
\end{eqnarray*}
where $R_\E$ (resp. $Ric_\E$) denotes the scalar curvature (resp. the Ricci curvature) of ($\E,g_\E)$ and $\mathcal{T}$ is the 
energy-momentum tensor which describes the matter content of the ambient spacetime. We also assume that $(\E^4,g_\E)$ satisfies the dominant energy condition that is its energy-momentum tensor
$\mathcal{T}$ has the property that, for every future directed causal
vector $\eta\in\Gamma(T\E)$, the vector field dual to the one form $-\mathcal{T}(\eta,.)$ is a future directed causal vector of $T\E$. \\


Let $M^3$ be an immersed spacelike hypersurface of $(\E^4,g_\E)$ with induced Riemannian metric $g$. Assume that $T$ is the future directed timelike normal vector 
to $M$ and denote by $K$ the associated second fundamental form defined by $K(X,Y)=g_\E(\nabla^\E_X T,Y)$, for all $X,Y\in\Gamma(TM)$. Here $\nabla^\E$ denotes the Levi-Civita connection of the spacetime. Then the Gauss, Codazzi and Einstein equations provide {\it constraint equations} 
on $M$ given by  
$$\left\lbrace
\begin{array}{rll}
\mu & = & \frac{1}{2}\big(R-|K|^2_M+({\rm Tr}_M(K))^2\big)\\
J & = & -\delta\big(K-{\rm Tr}_M(K)g\big) 
\end{array}
\right.
$$
where $R$ is the scalar curvature of $(M^3,g)$, $|K|^2$ and ${\rm Tr}(K)$ denote the squared norm and the trace of $K$ with respect to $g$ and $\delta$ is the divergence on $M$. Here $\mu$ (resp. $J$) 
is the energy (resp. the momentum) density of the matter fields given by 
\begin{eqnarray*}
\mu = \mathcal{T}(T,T)\quad\text{and}\quad J_i = \mathcal{T}(e_i,T)
\end{eqnarray*} 
for $1\leq i\leq 3$ and where $\{e_1,e_2,e_3\}$ is a local basis of the spatial tangent space of $M$. 
The dominant energy condition for the spacetime implies that $\mu\geq |J|$, as functions on $M$. A triplet $(M^3,g,K)$ which satisfies the dominant energy condition is called an {\it initial data set}.\\


Now we consider a codimension two spacelike orientable surface $\Sigma^2$ in the spacetime $\E^{4}$. We will represent by ${\mathcal H}$ the mean curvature vector field on $\Sigma^2$, defined as 
$$
{\mathcal H}={\rm tr\,}II,
$$
where $II$ is the second fundamental form of this immersion. Since the normal space at each point of $\Sigma^2$ is a Lorentzian
plane, it can be spanned by two future-directed null normal vector fields ${\mathcal N}_+$ and ${\mathcal N}_-$ normalized in such a way that $\langle{\mathcal N}_+,{\mathcal N}_-\rangle=-\frac1{2}$. We denote by $\theta_+$ and $\theta_-$ the components of ${\mathcal H}$ with respect to ${\mathcal N}_+$ and ${\mathcal N}_-$. They are the so-called future-directed null expansions of ${\mathcal H}$ and measure the area growth when $\Sigma^2$ varies in the corresponding directions. It is clear that
$$
|\mathcal{H}|^2= -\theta_+\theta_-.
$$ 
If $\theta_+$ and $\theta_-$ are both negative, the surface will be called a {\em trapped} surface. A surface with $\theta_+=0$ or $\theta_-=0$ is called  an {\em apparent horizon} (or a {\em marginally trapped} surface). Remark that if $\Sigma^2$ is trapped or marginally trapped, then the mean curvature vector ${\mathcal H}$ is a causal vector at each point. This is why, the fact that the mean curvature field ${\mathcal H}$ is spacelike everywhere is equivalent to that $\Sigma$ is an {\em untrapped} surface.  

In the case where $\Sigma^2$ spans a spacelike hypersurface in the spacetime, that is, when there exists a spacelike hypersurface
$\Omega^{3}$ immersed in $\E^{4}$ such that $\partial\Omega^{3}=\Sigma^2$, the normal null vector fields ${\mathcal N}_+$ and ${\mathcal N}_-$ may be ordered in such a way that they project onto directions tangent to $\Omega^{3}$ which are respectively {\em outer} and {\em inner} normal at each point of $\Sigma^2$. In other words, if $N$ is the inner normal unit vector field on $\Sigma^2$ tangent to $\Omega^{3}$ and $T$ is the future-directed timelike normal to $\Omega^{3}$ in $\E^{4}$, we put
$$
{\mathcal N}_+=\frac1{ 2}(T-N),\qquad {\mathcal N}_-=\frac1{2}(T+N).
$$
The second fundamental form of $\Sigma^2$ in $\E^{4}$ is given in terms of the Lorentzian basis of the normal bundle
provided by the hypersurface $\Omega^{}$ by
\begin{eqnarray*}
 II(X,Y)=g(AX,Y) N+g(BX,Y)T,
\end{eqnarray*}
for all $X,Y\in\Gamma(T\Sigma)$, where $AX:=-\nabla_XN$ denotes the shape operator of $\Sigma^2$ in $\Omega^{3}$ and $\nabla$ is the Levi-Civita connection of the Riemannian metric $g$ on $M$. The mean curvature vector field $\mathcal{H}$ of $\Sigma$ in $\E$ can be reexpressed by:
\begin{eqnarray*}
\mathcal{H}=\theta_+{\mathcal N}_-+\theta_-{\mathcal N}_+=H N+ {\rm tr}_\Sigma\,(K)\, T,
\end{eqnarray*}
where $H={\rm tr\,}A$ is the mean curvature of $\Sigma^2$ in $\Omega^{3}$ and ${\rm tr}_\Sigma\,(K)$ is the trace on $\Sigma^2$ of
the shape operator $K$ of $\Omega^{3}$ in ${\mathcal E}^{4}$. The norm of ${\mathcal H}$ can be also reexpressed: 
\begin{equation}\label{nullexp}
|\mathcal{H}|^2= H^2-{\rm tr}_\Sigma\,(K)^2=-\theta_+\theta_-,
\end{equation}
with $\theta_\pm={\rm tr}_\Sigma\,(K)\pm H$ are the future-directed null expansions of ${\mathcal H}$. The spacelike surfaces with $\theta_+<0$ (respectively, $\theta_-<0$) are referred to as {\em outer} (respectively, {\em inner}) {\em trapped} surfaces. It is easy to see that untrapped submanifolds, that is, codimension two spacelike submanifolds of a spacetime with spacelike mean curvature vector field, naturally divide into two disjoint classes.

\begin{lemma}\label{2classes}
Let $\Sigma^2$ be a compact spacelike codimension two submanifold embedded in a spacetime ${\mathcal E}^{4}$. Suppose that  its mean curvature vector field ${\mathcal H}$ is spacelike  and that $\Sigma^2$ is the boundary of a  spacelike  hypersurface  $\Omega^{3}$ in ${\mathcal E}^{4}$. Then $\Omega^{3}$ is either mean-convex or mean-concave.
\end{lemma}

\pf
It suffices to take into account that, if $(\theta_+,\theta_-)$ are the future-directed null expansions of the mean curvature vector field ${\mathcal H}$ associated to the embedding of $\Sigma^2$ in the domain $\Omega^{3}$, we have from (\ref{nullexp}) 
$$
0<|{\mathcal H}|^2=-\theta_+\theta_-,\qquad \theta_+-\theta_-=2H,
$$
where $H$ is the inner mean curvature function of $\Sigma^2$ in $\Omega^{3}$. The first of these two equalities implies that $\theta_+$ and $\theta_-$ have opposite signs everywhere on $\Sigma^2$. Then, from the second one, we have that either $H>0$ or $H<0$ on the whole 
of $\Sigma^2$.
\qed

Note that this fact obviously holds for higher dimensional initial data sets. In the following, an untrapped surface (resp. a codimension two untrapped submanifold) which bounds a compact, connected and mean-convex spacelike hypersurface will be referred to as an {\it outer untrapped surface} (resp. an {\it outer untrapped submanifold}). It is worth noting that round spheres in Euclidean slices are untrapped surfaces. The same occurs in general for large radial spheres in asymptotically flat spacelike hypersurfaces. 

We now give the precise statement of our main result:
\begin{theorem}\label{spacetime}
Let $\Omega$ be a compact domain with an outer untrapped boundary surface $\Sigma:=\partial\Omega$ in an initial data set $(M^3,g,K)$. If $\Omega$ has no apparent horizon in its interior, then for all $\varphi\in\Gamma(\sm\Sigma)$ we have:
\begin{eqnarray}\label{main-holo}
\int_\Sigma \Big(\frac{1}{|\mathcal{H}|}|\D \varphi|^2-\frac{|\mathcal{H}|}{4}|\varphi|^2\Big)\,d\Sigma\geq 0,
\end{eqnarray}
where $\sm\Sigma$ is the extrinsic spinor bundle on $\Sigma$ and $\D$ is the extrinsic Dirac operator (see Section $2$).
Moreover, if equality occurs then there exists a local isometric immersion of $\Omega$ in $\mathbb{R}^{3,1}$ with $K$ as second fundamental form.
\end{theorem}

As a direct application, we prove the following result
\begin{theorem}\label{liu-yau-general}
Under the same conditions of Theorem \ref{spacetime}, assume furthermore that $\Sigma$ admits an isometric and isospin immersion into $\mathbb{R}^3$ with mean curvature $H_0$. Then we have: 
\begin{eqnarray}\label{liu-yau-general-Ineq}
\int_{\Sigma}|\mathcal{H}|d\Sigma\leq\int_{\Sigma}\frac{H_0^2}{|\mathcal{H}|}d\Sigma.
\end{eqnarray}
Moreover, if equality occurs then $\Sigma$ is connected and there exists a local isometric immersion of $\Omega$ in $\mathbb{R}^{3,1}$ with second fundamental form given by $K$ and mean curvature vector of $\Sigma$ satisfying $|\mathcal{H}|=H_0$.  
\end{theorem}

If we consider the case of codimension two outer untrapped submanifolds in the $(n+2)$-dimensional Minkowski spacetime $\mathbb{R}^{n+1,1}$, we prove that we can remove the assumption on the non-existence of apparent horizons (see Theorem \ref{Mink}). Moreover in this situation, we completely characterize the equality case. Namely we have
\begin{theorem}\label{liu-yau-minkowski}
Let $\Sigma$ be a codimension two outer untrapped submanifold in $\mathbb{R}^{n+1,1}$. If $\Sigma$ admits an isometric and isospin immersion into $\mathbb{R}^{n+1}$ with mean curvature $H_0$, then Inequality  $(\ref{liu-yau-general-Ineq})$ holds and equality is achieved if and only if $\Sigma$ lies in a hyperplane in $\mathbb{R}^{n+1,1}$ and $\Sigma$ is connected.
\end{theorem}

\begin{remark}
{\rm In  Theorem \ref{liu-yau-general} and Theorem \ref{liu-yau-minkowski}, we assumed that the boundary hypersurface of a compact domain in a certain spin manifold admits an {\em isospin immersion} in a Euclidean space. In general, an $(n+1)$-dimensional spin manifold induces a spin structure on each of its orientable immersed hypersurfaces through their corresponding immersions (see Section \ref{hyper} below). Two distinct immersions of an orientable manifold $\Sigma^n$ into two (possibly different) $(n+1)$-dimensional spin manifolds are said to be {\em isospin} when the spin structures induced on $\Sigma^n$ from the corresponding ambient manifolds 
coincide (up to an equivalence). Recall that 
 spin structures on $\Sigma^n$ are parametrized by the cohomology group $H^1(\Sigma^n,{\mathbb Z}_2)$. Thus, for example, if $\Sigma^n$ is a simply-connected manifold, any two immersions of $\Sigma^n$ in two arbitrary $(n+1)$-dimensional spin manifolds must be isospin. Consequently, if the surface $\Sigma$ in Theorem \ref{liu-yau-general} has genus zero or the hypersurface $\Sigma$ in Theorem \ref{liu-yau-minkowski} is simply-connected, we only need to suppose that they are mean-convex in their initial data sets and that they can be immersed as hypersurfaces in a Euclidean space. 
 
 Also, it is clear that, when the two immersions defined on $\Sigma^n$ lie in a same ambient space and are {\em regularly homotopic}, the associated induced spin structures  are equivalent. In fact,  two immersions are said to be regularly homotopic ({\em isotopic}, according to Pinkall and others, see \cite{Pi}) if we may pass continuously from one to the other through a family of immersions. Consequently, they determine the same class in $H^1(\Sigma^n,{\mathbb Z}_2)$. Indeed, in the case $n=2$, two spin structures induced from the spin structure of the $3$-dimensional spin ambient space through two different embeddings are equivalent if and only if they are regularly homotopic (besides of \cite{Pi}, see \cite[pp. 104--105]{HH} and \cite[p. 656]{BeS}).

Then, take any compact mean-convex surface $\Sigma$  embedded in ${\mathbb R}^3$. This surface bounds a compact domain in the three-dimensional Euclidean space which is a totally geodesic initial data set in the Minkowski space ${\mathbb R}^{3,1}$. If we slightly deform this surface the positivity of the mean curvature is preserved by continuity and, from the arguments above, the same holds for the induced spin structure. So, there are examples of mean-convex boundaries in initial data sets of spacetimes admitting isospin immersions in Euclidean spaces. Many of them are non-convex. In fact, take $\Sigma$ to be, for instance, a right cylinder with two half-spheres closing its extremes (after smoothing) or a torus of revolution thin enough (if we want to have some point with negative Gau{\ss} curvature).   

Note that, if $\Sigma$ is not convex, we cannot use the Weyl theorem and so we do not know whether it is possible to immerse $\Sigma$ isometrically in Euclidean space ${\mathbb R}^3$ or not. This is why, in this case, Theorem \ref{liu-yau-general} and Theorem \ref{liu-yau-minkowski} should be viewed as a comparison theorems for the mean curvatures of two immersions in the spirit of a classical result by Herglotz. Indeed, in $1934$, Herglotz \cite{He} gave a succinct proof of Cohn-Vossen's rigidity result for convex surfaces based on an integral inequality involving the second fundamental forms of two embedding  (see, for example, \cite[Section 7.4]{MR}). Our Theorem \ref{liu-yau-general} provides an inequality of this type which could be a first step in order to enlarge the Cohn-Vossen theorem to include Euclidean mean-convex compact surfaces.

In this direction, one can easily see that Theorem 3  implies that the integral of the mean curvature is preserved through {\em bendings} of compact mean-convex hypersurfaces embedded in a Euclidean space.  This was first proved by Almgren and Rivin (\cite{AR}, see
also \cite{RS}).}\\
\end{remark}

Recall that, in \cite{LY2} (see also \cite{LY1}), Liu and Yau proved the following positivity result:
{\emph {Let $(\Omega^3,g,K)$ be an initial data set for the Einstein equation. Suppose that the boundary $\partial\Omega$ has finitely many components $\Sigma_i$, $1\leq i\leq l$, each of which has positive Gauss curvature and spacelike mean curvature vector in the spacetime. Then for all $i$:
\begin{eqnarray}\label{liu-yau-ineq}
\int_{\Sigma_i}|\mathcal{H}|d\Sigma\leq\int_{\Sigma_i} H_0d\Sigma.
\end{eqnarray}
Moreover if equality occurs for some $i\in\{1,...,l\}$, then $\partial\Omega$ is connected and the spacetime is flat along $\Omega$.}}\\

The proof of this result relies on a generalized version of the Positive Mass Theorem and on the resolution of the Jang equation. One of the key ingredients in their proof is provided by the Weyl embedding theorem \cite{We} which asserts that $\Sigma$ embeds isometrically as a strictly convex hypersurface in $\mathbb{R}^3$ is equivalent to the fact that $\Sigma$ has positive Gauss curvature. 
Note that by the Cauchy-Schwarz inequality, Inequality (\ref{liu-yau-ineq}) implies  (\ref{liu-yau-general-Ineq}).\\

More recently, Eichmair, Miao and Wang \cite{EMW} generalized Inequality (\ref{liu-yau-ineq}) for time-symmetric initial data under weaker convexity assumptions for the embedding of $\Sigma$ in $\mathbb{R}^3$. We point out that, in contrast to Liu-Yau's result, we do not assume that the immersion is a {\it strictly convex embedding}. In particular, the mean curvature $H_0$ is not assumed to be positive.


\section{The Riemannian setting}



\subsection{Preliminaries on Spin Manifolds} 


Let $(M,g)$ be an $(n+1)$-dimen\-sional Riemannian Spin manifold, which we will suppose from now on to be connected, and denote by  ${\nabla}$ the Levi-Civita connection on its tangent bundle $TM$. We choose a Spin structure on $M$ and consider the corresponding spinor bundle $\SM$, a rank $2^{\left[\frac{n+1}{2}\right]}$ complex vector bundle. Denote by $\gamma$ the Clifford multiplication
\begin{equation}\label{Clm}
\gamma:\C\ell(M)\longrightarrow \hbox{End}(\SM)
\end{equation}
which is a fiber preserving algebra morphism. Then $\SM$ becomes a bundle of complex left modules over the Clifford bundle $\C\ell(M)$ over the manifold $M$. When $(n+1)$ is even, the spinor bundle splits into the direct sum of the {\em positive} and {\em negative} chiral subbundles
\begin{equation}\label{chiral}
\SM=\SM^+\oplus\SM^-,
\end{equation}
where $\SM ^{\pm}$ are defined to be the $\pm 1$-eigenspaces of the endomorphism $\gamma(\omega_{n+1})$, with $\omega_{n+1}=i^{\left[\frac{n+2}{2}\right]}e_1\cdot e_2\cdots e_{n+1}$, the complex volume form. 

On the spinor bundle $\SM$, one has (see \cite{LM}) a natural Hermitian metric, denoted by $\langle\; ,\;\rangle$, and the spinorial Levi-Civita  connection ${\nabla}$ acting on spinor fields. It is well-known that the Hermitian scalar product, the Levi-Civita connection ${\nabla}$ and the Clifford multiplication (\ref{Clm}) satisfy, for any spinor fields $\psi,\varphi\in\Gamma(\SM)$ and any 
tangent vector fields $X,Y\in \Gamma (TM)$, the following compatibility conditions: 
\begin{eqnarray}
&\langle\gamma(X)\psi,\gamma(X)\varphi\rangle =|X|^2
\langle\psi,\varphi\rangle&\label{comp2}\\
& X\langle\psi,\varphi\rangle = \langle{\nabla}_X\psi,
\varphi\rangle+\langle\psi,{\nabla}_X\varphi
\rangle&\label{comp1}\\  
&{\nabla}_X\big(\gamma(Y)\psi\big) 
= \gamma({\nabla}_XY)\psi
+\gamma(Y){\nabla}_X\psi.&\label{comp3}
\end{eqnarray}
Since ${\nabla} \omega_{n+1}=0$, for $(n+1)$ even, the decomposition (\ref{chiral}) is orthogonal and ${\nabla}$ preserves this decomposition.

The Dirac operator ${D}$ on $\SM$ is the first order elliptic differential operator locally given by 
$$
{D}=\sum_{i=1}^{n+1}\gamma(e_i){\nabla}_{e_i},
$$
where $\{e_1,\dots,e_{n+1}\}$ is a local orthonormal frame of $TM$. When  $(n+1)$ is even, the Dirac operator interchanges positive and negative spinor fields, that is,
\begin{eqnarray*}
{D} : \Gamma(\SM^\pm)\longmapsto \Gamma(\SM^\mp).
\end{eqnarray*}


\subsection{Hypersurfaces and induced Structures}\label{hyper}


In this section, we compare the restriction $\sm \Sigma$ of the spinor bundle $\S M$ of a Spin manifold $M$ to an orientable 
hypersurface $\Sigma$ immersed into $M$ and its Dirac-type operator $\D$ to the intrinsic spinor bundle $\S\Sigma$ of the induced Spin structure on $\Sigma$ and its fundamental Dirac operator $D_\Sigma$. A fundamental case will be when the hypersurface $\Sigma$ is just the boundary $\partial M$ of a manifold $M$. These facts are in general well-known (see for example \cite{Bu,Tr, Ba2,BFGK,HMZ1,HMZ2,HMZ3,HM1}). For completeness, we introduce the notations and the key facts.

Denote by $\nb$ the Levi-Civita connection associated with the induced Riemannian metric on $\Sigma$. The Gau{\ss} formula says that
\begin{equation}\label{riem-gaus}
\nb_XY={\nabla}_XY-g(AX,Y) N,
\end{equation}
where $X,Y$ are vector fields tangent to the hypersurface $\Sigma$, the vector field $N$ is a global unit field normal to $\Sigma$ and $A$ stands for the shape operator corresponding to $N$, that is,
\begin{equation}\label{shap-oper}
{\nabla}_XN=-A X,\qquad \forall X\in \Gamma(T\Sigma).
\end{equation}
We have that the restriction
\begin{equation}
\sm\Sigma:=\SM_{|\Sigma}\nonumber
\end{equation}
is a left module over $\C\ell(\Sigma)$ for the induced Clifford multiplication$$
\mult:\C\ell(\Sigma)\longrightarrow \hbox{End}(\sm\Sigma)$$
given by 
\begin{equation}\label{Clmind}
\mult(X)\psi=\gamma(X)\gamma(N)\psi
\end{equation}
for every $\psi\in\Gamma(\sm\Sigma)$ and $X\in\Gamma(T\Sigma)$ (note that a spinor field on the ambient manifold $M$ and its restriction to the hypersurface $\Sigma$ will be denoted by the same symbol). Consider on $\sm\Sigma$ the Hermitian metric $\langle\; ,\;\rangle$ induced from that of $\SM$. This metric immediately satisfies the compatibility condition (\ref{comp2}) if one considers on $\Sigma$
the Riemannian metric induced from $M$ and the Clifford multiplication $\mult$  defined in (\ref{Clmind}). Now the Gauss formula 
(\ref{riem-gaus}) implies that the Spin connection $\nb$ on $\sm\Sigma$ is given by the following spinorial Gauss formula
\begin{equation}\label{spin-gaus}
\nb_X\psi={\nabla}_X\psi-\frac{1}{2}\mult(A X)\psi
={\nabla}_X\psi-
\frac{1}{2}\gamma(A X)\gamma(N)\psi\,
\end{equation}
for every $\psi\in\Gamma(\sm\Sigma)$ and $X\in\Gamma(T\Sigma)$. Note that the compatibility conditions (\ref{comp2}), (\ref{comp1}) and (\ref{comp3}) are satisfied by $(\sm\Sigma,\mult,\langle\;,\;\rangle,\nb)$. 

Denote by ${\D}:\Gamma(\sm\Sigma)\rightarrow \Gamma(\sm\Sigma)$ the Dirac operator associated with the Dirac bundle $\sm\Sigma$ over the hypersurface. It is a well-known fact that ${\D}$ is a first order elliptic differential operator which is formally $L^2$-selfadjoint. By (\ref{spin-gaus}), for any spinor field $\psi\in\Gamma(\SM)$, we have 
\begin{equation}
{\D}\psi=\sum_{j=1}^n\mult(e_j)\nb_{e_j}\psi
=\frac{1}{2}H\psi-\gamma(N)\sum_{j=1}^n\gamma(e_j)
{\nabla}_{e_j}\psi,\nonumber
\end{equation}
where $\{e_1,\dots,e_n\}$ is a local orthonormal frame of $T\Sigma$ and $H=\hbox{trace\,}A$ is the mean curvature of $\Sigma$ corresponding to the orientation $N$. Using (\ref{spin-gaus}) and (\ref{shap-oper}), it is straightforward to see that the skew-commu\-ta\-ti\-vity rule
\begin{equation}\label{D-commutes}
\D\big(\gamma(N)\psi\big)=-\gamma(N)\D\psi
\end{equation}
holds for any spinor field $\psi\in\Gamma({\sm}\Sigma)$. It is important to point out that, from this fact, {\em the spectrum of $\D$ is always symmetric with respect to zero}, while this is the case for the Dirac operator $D_\Sigma$ of the intrinsic spinor bundle {\em only when $n$ is even}. Indeed, in this case, we have an isomorphism of Dirac bundles
$$
({\sm} \Sigma,\mult,\D)\equiv (\bS \Sigma,\gamma_\Sigma,D_\Sigma)
$$
and the decomposition ${\sm} \Sigma={\sm}\Sigma^+\oplus {\sm} \Sigma^-$, given by ${\sm} \Sigma^\pm:=\{\psi\in{\sm} \Sigma\,|\,i\gamma(N)\psi=\pm\psi\},$ corresponds to the chiral decomposition of the spinor bundle $\bS \Sigma$. Hence $\D$ interchanges 
${\sm}\Sigma^+$ and ${\sm}\Sigma^-$.

When $n$ is odd the spectrum of $D_\Sigma$ is not necessarily  symmetric. In fact, in this case, the spectrum of $\D$ is just the symmetrization of the spectrum of $D_\Sigma$. This is why the decomposition of $\S{M}$ into positive and negative chiral spinors induces an orthogonal and $\mult,\D$-invariant decomposition ${\sm}\Sigma={\sm}\Sigma_+\oplus {\sm}\Sigma_-$, with ${\sm}\Sigma_\pm:=(\S{M}^\pm)_{|\Sigma}$, in such a way that
$$
({\sm}\Sigma_\pm,\mult,\D_{|{\sm}\Sigma_\pm}) \equiv (\bS \Sigma, \pm\gamma_\Sigma, \pm D_\Sigma).
$$ 
Moreover, $\gamma(N)$ interchanges the decomposition and both maps $\gamma(N) : {\sm}\Sigma_\pm \longrightarrow {\sm}\Sigma_\mp$
are isomorphisms.

Consequently, to study the spectrum of the induced operator $\D$ is equivalent to study the spectrum of the Dirac operator $D_\Sigma$ of the Riemannian Spin structure induced on the hypersurface $\Sigma$.


\subsection{A spinorial Reilly type Inequality for Manifolds with Boundary}


In this section, we prove a spinorial Reilly type inequality (\cite{LY1},\cite{Ra}). 

Recall that, on a compact $(n+1)$-dimensional Riemannian Spin manifold $M$ with boundary $\Sigma=\partial M$, for any spinor field $\psi\in\Gamma(\bS M)$, the fundamental Schr{\"o}dinger-Lichnerowicz formula is given by:
\begin{equation}\label{boun-weit-twis-ineq}
\int_\Sigma\big(\langle {\D}\psi,\psi\rangle-\frac{H}{2}|\psi|^2\big)\,d\Sigma= 
\int_ M \left(\frac{1}{4} {R}|\psi|^2 + |\nabla\psi|^2 - |{D}\psi|^2\right)\,d M
\end{equation}
where $R$ is the scalar curvature of $M$. Note that, the assumption $R\ge 0$ is quite natural and has been used intensively to get, in particular, lower bounds on both $D$ and ${\D}$.  However, in our situation (see Section \ref{jang-equation}), we have a weaker assumption on the scalar curvature. More precisely, we assume that there exits a smooth vector field $X\in\Gamma(TM)$ such that:
\begin{eqnarray}\label{scalarcond}
R\geq 2|X|^2+2\delta(X)
\end{eqnarray}
where $|X|^2=g(X,X)$ and $\delta$ is the divergence of $X=\sum_{j=1}^nX^je_j\in\Gamma(TM)$, locally given by
\begin{eqnarray*}
\delta(X)=-\sum_{i=1}^{n+1}{e_i} (X^i).
\end{eqnarray*} 
Then we prove an adapted Reilly type inequality. Namely:
\begin{proposition}\label{modlich}
Let $M$ a compact Riemannian Spin manifold with boundary $\Sigma$ such that there exists a smooth vector field $X\in\Gamma(TM)$ satisfying 
(\ref{scalarcond}), then 
\begin{eqnarray}\label{mod-weit-ineq}
\int_{\Sigma}\<\D\psi-\frac{1}{2}\big(H+g(X,N)\big)\psi,\psi\>\,d\Sigma\geq\int_{M}\big(\frac{1}{2}|\nabla\psi|^2-|D\psi|^2\big)\,dM.
\end{eqnarray}
Moreover equality occurs if and only if the spinor field $\psi$ satisfies 
\begin{eqnarray}\label{reilly-equa}
\nabla_Y\psi=-g(X,Y)\psi
\end{eqnarray}
for all $Y\in\Gamma(TM)$.
\end{proposition}
\pf
First note that, since 
\begin{eqnarray*}
 \delta(|\psi|^2 X)=-X(|\psi|^2)+|\psi|^2\delta(X),
\end{eqnarray*}
Stokes formula gives
\begin{eqnarray*}
\int_M\frac{R}{4}|\psi|^2 dM& = & \int_M\big(\frac{R}{4}-\frac{1}{2}\delta(X)\big)|\psi|^2dM+\frac{1}{2}\int_M\delta(X)|\psi|^2dM\\
& = & \frac{1}{4}\int_M\big(R-2\delta(X)\big)|\psi|^2dM+\frac{1}{2}\int_M X(|\psi|^2)dM\\
& & +\frac{1}{2}\int_\Sigma g(X,N)|\psi|^2d\Sigma.
\end{eqnarray*}
Inserting this identity in (\ref{boun-weit-twis-ineq}) leads to:
\begin{eqnarray*}
&&\int_{\Sigma}\<\D\psi- \frac{1}{2}\big(H + g(X,N)  \big) , \psi\>d\Sigma  =  \\
&& \int_M\big(\frac{1}{4}\big(R-2\delta(X)\big)|\psi|^2+\frac{1}{2}X(|\psi|^2)\big)dM
 +\int_{M}\big(|\nabla\psi|^2-|D\psi|^2\big)dM
\end{eqnarray*}
and using (\ref{scalarcond}), we conclude that:
\begin{eqnarray}\label{schrint}
\int_{\Sigma}\<\D\psi-\frac{1}{2}\big(H+g(X,N)\big) , \psi\>d\Sigma & \geq &
\int_{M}\big(\frac{1}{2}|X|^2|\psi|^2+\frac{1}{2}X(|\psi|^2)\big)dM\nonumber\\
& & +\int_{M}\big(|\nabla\psi|^2-|D\psi|^2\big)dM.
\end{eqnarray}
If we let $\widetilde{\nabla}_Y\psi:=\nabla_Y\psi+g(X,Y)\psi$, it is straightforward to compute
\begin{eqnarray*}
|\widetilde{\nabla}\psi|^2=|\nabla\psi|^2+|X|^2|\psi|^2+2{\rm Re}\<\nabla_X\psi,\psi\> 
\end{eqnarray*}
and since $2{\rm Re}\<\nabla_X\psi,\psi\>=X(|\psi|^2)$, we get:
\begin{eqnarray*}
\frac{1}{2}X(|\psi|^2)\geq -\frac{1}{2}|\nabla\psi|^2-\frac{1}{2}|X|^2|\psi|^2
\end{eqnarray*}
with equality if and only if $\widetilde{\nabla} \psi =0$. Combining this last inequality with (\ref{schrint}) finishes the proof.
\qed


\subsection{A Local Boundary Elliptic Condition for the Dirac Operator}\label{sec:gcbc}


As before, $\Sigma$ is the boundary of an $(n+1)$-dimensional Riemannian Spin compact manifold $M$. We define
two pointwise projections 
$$
P_\pm:\sm\Sigma\longrightarrow \sm\Sigma
$$
on the induced Dirac bundle over the hypersurface, as follows
\begin{equation}\label{defP}
P_\pm=\frac{1}{2}\big(\hbox{Id}_{\sm\Sigma}\pm i\gamma(N)\big).
\end{equation}

It is a well known fact that these two orthogonal projections $P_\pm$ acting on the spin bundle $\sm\Sigma$ provide local elliptic boundary conditions for the Dirac operator ${D}$ of $M$. The ellipticity of these boundary conditions and that of the Dirac operator $D$, allow to solve boundary-value problems for $D$ on $M$ by prescribing on the boundary $\Sigma$, the corresponding $P_\pm$-projections of the solutions. Namely we have:
\begin{proposition}\cite{HM1}\label{boun-prob2}
Let $M$ be a compact Riemannian Spin manifold with boundary, a hypersurface $\Sigma$. If $\varphi\in\Gamma(\sm\Sigma)$ is a smooth spinor field of the induced Dirac bundle, then the following boundary-value problem for the Dirac operator 
$$ 
\left\{
\begin{array}{lll}
{D}\psi&=0 \qquad&\hbox{ {\rm on} } M   \\
P_\pm(\psi_{|\Sigma})&=P_\pm\varphi \qquad&\hbox{ \rm on }\Sigma
\end{array}
\right. 
$$
has a unique smooth solution $\psi\in\Gamma(\bS M)$.
\end{proposition}

For a more general discussion on boundary conditions for the Dirac operator, we refer to \cite{BW}, \cite{BB} or \cite{BC}.


\subsection{A Holographic Principle for the Existence of Parallel Spinors}


It is by now standard (see \cite{HMZ2, HMZ3}) to make use of (\ref{boun-weit-twis-ineq}) for a compact Riemannian Spin manifold $M$
with non-negative scalar curvature $R$, together with the solution of an appropriate boundary-value problem for the  Dirac operator $D$ of $M$, in order to establish a certain integral inequality for the induced Dirac operator $\D$ of the boundary hypersurface $\partial M=\Sigma$. In \cite{Ra}, the third author uses such arguments for compact manifolds whose scalar curvature satisfies (\ref{scalarcond}). In this section, we generalize the holographic principle for the existence of parallel spinors proved by the first two authors in \cite{HM1} in the context studied in \cite{Ra}. 

First, we need to recall the following fact:
\begin{lemma}\cite{HMZ3}\label{lemma1} 
For any smooth spinor field $\psi\in\Gamma(\sm\Sigma)$ we have
$$
\int_\Sigma\langle{\D}\psi,\psi\rangle\,d\Sigma=2 \int_\Sigma
\langle{\D}P_+\psi,P_-\psi\rangle\,d\Sigma.
$$
\end{lemma}
The proof simply relies on the self-adjointness of the Dirac operator $\D$ and on the identities
\begin{equation}\label{antiDP}
{\D}P_\pm=P_\mp{\D}
\end{equation}
which are obtained using (\ref{D-commutes}) and (\ref{defP}).

\begin{proposition}\label{proposition4} 
Let $M$ be a compact Riemannian  Spin manifold with scalar curvature satisfying $(\ref{scalarcond})$ and such that
\begin{eqnarray*}
F:=H+g(X,N)>0.
\end{eqnarray*}
For any  $\varphi\in\Gamma(\sm\Sigma)$, one has 
\begin{equation}\label{ineq+}
0\le \int_\Sigma \big(\frac{1}{F}|\D P_+\varphi|^2-\frac{F}{4}|P_+\varphi|^2\big)\,d\Sigma.
\end{equation}
Moreover, equality holds if and only if there exists a parallel spinor field $\psi
\in \Gamma({\mathbb S} M)$ such that $P_+\psi=P_+\varphi$ along the boundary hypersurface $\Sigma$ and the vector field $X$ vanishes identically on $M$.
\end{proposition}
\pf
Take any spinor field $\varphi\in\Gamma(\sm\Sigma)$ of the induced spinor bundle on the hypersurface and consider the following boundary-value problem
$$
\left\{
\begin{array}{lll}
{D}\psi&=0 \qquad&\hbox{ {\rm on} } M \\
P_+\psi&=P_+ \varphi \qquad&\hbox{ {\rm on} }\Sigma 
\end{array}
\right.
$$ 
for the Dirac operator ${D}$ and the boundary condition $P_+$. The existence and uniqueness of a smooth solution $\psi\in \Gamma(\bS M)$ for this boundary problem is ensured by Proposition \ref{boun-prob2}. This solution $\psi$ inserted in Inequality (\ref{mod-weit-ineq}), translates to
\begin{equation}\label{final-ineq}
0\leq\frac{1}{2}\int_{M}|\nabla\psi|^2dM\le \int_\Sigma\big(\langle {\D}\psi,\psi\rangle
-\frac{F}{2}|\psi|^2\big)\,d\Sigma.
\end{equation}
Note that if equality is achieved, then $\psi$ is a parallel spinor field satisfying (\ref{reilly-equa}). Since such a spinor field has no zeros, the vector field $X$ vanishes identically on the whole of $M$. Inequality (\ref{final-ineq}) combined with Lemma \ref{lemma1} together with the fact that the decomposition 
$$\psi=P_+\psi+P_-\psi$$
is pointwise orthogonal, imply
\begin{equation}\label{final-ineq-P}
0\le \int_\Sigma \left(2\langle{\D}P_+\psi,P_-\psi\rangle
-\frac{F}{2}|P_+\psi|^2-\frac{F}{2}|P_-\psi|^2\right)\,d\Sigma.
\end{equation}
The function $F$ being assumed positive on $\Sigma$, it follows
\begin{eqnarray}
&{\displaystyle
0\le \big| \sqrt{\frac{2}{F}}\D P_+\psi-\sqrt{\frac{F}{2}} P_-\psi\big|^2
= } \nonumber \\
&{\displaystyle \frac{2}{F}|\D P_+\psi|^2+\frac{F}{2}|P_-\psi|^2-2\langle\D P_+\psi,
P_-\psi\rangle.}\nonumber
\end{eqnarray}
In other words, we have
\beQ
2\langle\D P_+\psi , P_-\psi\rangle -\frac{F}{2}|P_-\psi|^2 \le \frac{2}{F}|\D P_+\psi|^2,
\eeQ
which, when combined with Inequality (\ref{final-ineq-P}), implies Inequality (\ref{ineq+}). Now if equality holds, we already noticed that the spinor field $\psi$ must be parallel with $P_+\psi=P_+\varphi$ and $X\equiv 0$. 

Conversely, if we assume that there is a parallel spinor field $\psi$ on $M$ and $X\equiv 0$ then we are in the situation covered in \cite{HM1}. 
\qed

With this, we are ready to state the main result of this section: 
\begin{theorem}\label{holo-princ}
Let $M$ be a compact Riemannian Spin $(n+1)$-dimen\-si\-o\-nal manifold and $X\in\Gamma (TM)$ such that 
\begin{eqnarray*}
R\geq 2|X|^2+2\delta(X) \quad\text{and}\quad F:=H+g(X,N)>0.
\end{eqnarray*}
Then, for any spinor field $\varphi\in\Gamma(\sm\Sigma)$, one has
\begin{equation}\label{ineq-total}
0\le \int_\Sigma \big(\frac{1}{F}|\D \varphi|^2-\frac{F}{4}|\varphi|^2\big)\,d\Sigma.
\end{equation}
Equality holds if and only if there exist two parallel spinor fields $\Psi^+,\Psi^-\in\Gamma({\mathbb S} M)$ such that $P_+\Psi^+=P_+\varphi$ and $P_-\Psi^-=P_-\varphi$
on the boundary and $X\equiv 0$.
\end{theorem}
\pf
From the symmetry between the two boundary conditions $P_+$ and $P_-$ for the Dirac operator on $M$ (see Proposition \ref{boun-prob2} and Lemma \ref{lemma1}), one can repeat the proof of Proposition \ref{proposition4} to get the inequality corresponding  to (\ref{ineq+}) where the {\em positive} projection $P_+$ is replaced by the {\em negative} one $P_-$. Hence, for any spinor field $\varphi\in\Gamma(\sm\Sigma)$, we also have
\begin{equation}\label{ineq-}
0\le \int_\Sigma \big(\frac{1}{F}|\D P_-\varphi|^2-\frac{F}{4}|P_-\varphi|^2\big)\,d\Sigma.
\end{equation}
Taking into account the relation (\ref{antiDP}) and the pointwise orthogonality of the projections $P_\pm$, the sum of the two inequalities (\ref{ineq+}) and (\ref{ineq-}) yields (\ref{ineq-total}). The equality case is a consequence of Proposition \ref{proposition4}.
\qed

\begin{remark}
{\rm Note that, as observed in \cite{HM1}, equality in (\ref{ineq-total}), does not imply that the two parallel spinors in Theorem \ref{holo-princ} coincide.} 
\end{remark}

\noindent We should also mention that Inequality (\ref{ineq-total}) has a nice interpretation in terms of the first eigenvalue of the boundary Dirac operator $\D_F$ associated with the conformal metric $g_F=F^2g$. More precisely, we have:
\begin{corollary}
Let $(M^{n+1},g)$ be an $(n+1)$-dimensional compact connected Riemannian Spin manifold satisfying the assumptions of Theorem \ref{holo-princ}. Then, the first non-negative eigenvalue $\lambda_1(\D_F)$, of the Dirac operator corresponding to the conformal metric $g_F=F^2g$, satisfies
\begin{equation}
\lambda_1(\D_F)\ge \frac{1}{2}\nonumber
\end{equation}
and equality holds if and only if $M$ admits a non trivial parallel spinor (and  $X\equiv 0$). In this case, the  eigenspace corresponding to $\lambda_1(\D_{F})=\frac{1}{2}$ consists of  restrictions to $\Sigma$ of parallel spinor fields on
$M$ multiplied by the function $F^{-\frac{n-1}{2}}$. Furthermore, the boundary hypersurface $\Sigma$ has to be connected.
\end{corollary}

The proof is omitted since it is similar to that of Theorem $1$ in \cite{HM1}.


\subsection{A Discussion on Quasi-Local Masses}


In this section, we consider a $3$-dimensional compact  connected Riemannian manifold $(M^3,g)$ with non-negative scalar curvature whose boundary $\Sigma^2$ has positive mean curvature $H$. Note that, since  $M$ is a $3$-dimensional manifold, it is necessary Spin.  Moreover, we also assume that there exits an immersion $\iota_0$ of the surface $\Sigma$ in $\mathbb{R}^3$ with mean curvature $H_0$. \\

One of the fundamental  results in classical general relativity is certainly the proof of the positivity of the total energy by Schoen-Yau \cite{sy1} and Witten \cite{Wi}. This led to the more ambitious claim to associate energy to extended, but finite, spacetime domains, i.e., at the quasi-local level. Obviously, the quasi-local data could provide a more detailed characterization of the states of the gravitational field than the global ones, so they are interesting for their own right. For a complete review of these topics, we refer to \cite{Sz}. It is currently required that a quasi-local mass satisfies natural properties, among which:
\begin{enumerate}[(I)]
\item\label{nonnegativity} {\bf Non-negativity}: $\mg\geq 0$;
\item\label{rigidity} {\bf Rigidity}: $\mg=0$ if and only if $\Sigma$ is in the Minkowski spacetime;
\item\label{monotonicity} {\bf Monotonicity}: If $\Sigma_1=\partial M_1$ and $\Sigma_2=\partial M_2$ with $M_1\subset M_2$, then $\mathcal{M}(\Sigma_1)\leq\mathcal{M}(\Sigma_2)$;
\item\label{ADMlimit} {\bf ADM limit}: If $(\Sigma_k)$ is a sequence of surfaces that exhaust an asymptotically flat manifold $(N^3,g)$ then 
\begin{eqnarray*}
\lim_{k\rightarrow\infty}\mathcal{M}(\Sigma_k)=m_{ADM}(g)
\end{eqnarray*}

where $m_{ADM}(g)$ is the ADM mass of $(N,g)$.
\item\label{blackhole} {\bf Black hole limit}: If $\Sigma$ is a horizon in an asymptotically flat manifold $(N^3,g)$, then 
\begin{eqnarray*}
\mg=\sqrt{\frac{A}{16\pi}}
\end{eqnarray*}
where $A$ is the area of $\Sigma$.
\end{enumerate}

In \cite{BY}, Brown and York proposed the following definition for the quasi-local mass of a surface $\Sigma$ (now called the Brown-York mass):
\begin{eqnarray*}
\mby:=\frac{1}{8\pi}\int_{\Sigma}(H_0-H)d\Sigma.
\end{eqnarray*}
The non-negativity of $\mby$ is proved in \cite{ST1} under additional assumptions. Indeed, they impose that $\iota_0$ is a {\it strictly convex} isometric embedding which by the Weyl embedding theorem \cite{We} is equivalent to the fact that $\Sigma$ has positive Gauss curvature. Moreover in this situation, the embedding $\iota_0$ is unique up to an isometry of $\mathbb{R}^3$. \\
 
Recently, Lam \cite{L} proposed in his thesis the following definition:
\begin{eqnarray*}
\ml:=\frac{1}{16\pi}\int_{\Sigma}\frac{1}{H_0}(H^2_0-H^2)d\Sigma.
\end{eqnarray*}
He proves that $\ml$ has several interesting properties for certain surfaces in complete asymptotically flat Riemannian manifolds that are the graphs of smooth functions over $\mathbb{R}^3$ (see \cite{L} for a precise description). More precisely, it satisfies Properties (\ref{nonnegativity}), (\ref{monotonicity}), (\ref{ADMlimit}) and (\ref{blackhole}). Moreover, using the Cauchy-Schwarz inequality it is straightforward to check that $\mby\geq\ml$.

From the work of the first two authors \cite{HM1}, we can define a quasi-local mass similar to Brown-York and Lam and prove its non-negativity in the more general context described in the beginning of this section. Indeed, if we let
\begin{eqnarray*}
\mhm:=\frac{1}{16\pi}\int_{\Sigma}\frac{1}{H}(H^2_0-H^2)d\Sigma
\end{eqnarray*}
then from the immersion $\iota_0$, there exists a spinor field $\Psi_0\in\Gamma(\sm\Sigma)$ satisfying the following Dirac equation
\begin{eqnarray*}
\D\Psi_0=\frac{H_0}{2}\Psi_0,\qquad |\Psi_0|=1.
\end{eqnarray*}
It is obtained by taking the restriction to $\Sigma$ of a parallel spinor field on $\mathbb{R}^3$. Now taking $\Psi_0$ in Inequality (\ref{ineq-total}) with $X\equiv 0$ and $F=H$ gives $\mhm\geq 0$. Moreover, from \cite{HM1}, $\mhm=0$ if and only if $M$ is a Euclidean domain and the embedding of $\Sigma$ in $M$ and its immersion in $\mathbb{R}^{3}$ are congruent. In other words, Properties (\ref{nonnegativity}) and (\ref{rigidity}) are satisfied. \\

Note that if we assume that $\Sigma$ has positive Gauss curvature (which is a stronger assumption) then using the Cauchy-Schwarz inequality implies that $\mhm\geq\mby$ and the nonnegativity of $\mhm$ follows from the nonnegativity of the Brown-York mass. On the other hand, it is also proved in \cite{HM1} (see the proof of Corollary $10$) that (\ref{ADMlimit}) holds. However it is clear from the definition that the mass $\mhm$ is not defined for minimal surfaces (and so for apparent horizons). Moreover, the monotonicity Property (\ref{monotonicity}) is not satisfied in general. Take for example the $3$-dimensional Schwarzschild manifold $(N^3,g)=(\mathbb{R}^3\setminus\{0\},u^4g_{eucl})$ where $u:=1+\frac{M}{2r}$, $M>0$ and $g_{eucl}$ is the Euclidean metric. For a sphere $\mathbb{S}^2_r$ in $N^3$, its isometric image in $\mathbb{R}^3$ is $\mathbb{S}^2_{ru^2}$. Thus $H_0=\frac{2}{ru^2}$ and since the Schwarzschild metric is conformal to the Euclidean metric, we have
\begin{eqnarray*}
H=u^{-2}\big(\frac{2}{r}+\frac{4}{u}\frac{\partial u}{\partial r}\big).
\end{eqnarray*} 
A direct computation gives
\begin{eqnarray*}
m(\mathbb{S}^2_r)=M\frac{r+\frac{M}{2}}{r-\frac{M}{2}}
\end{eqnarray*}
and so $m(\mathbb{S}^2_r)$ is monotonically {\it decreasing} to the $ADM$ mass $M$ as $r$ goes to infinity.


\section{Spacelike surfaces in initial data sets}



\subsection{The Jang Equation}\label{jang-equation}


In this section, we recall some well-known facts on the Jang equation (for more details, we refer to \cite{sy1}, \cite{yau1} or \cite{AEM}). This equation first appears in \cite{jang} in his attempt to prove the positive mass theorem using the inverse mean curvature flow. However, as shown by Schoen and Yau \cite{sy1}, this equation can be used to reduce the proof of the general positive mass theorem to the case of time-symmetric initial data sets (that is $K_{ij}=0)$ previously obtained by the same authors in \cite{symac}. More recently,  Liu and Yau (\cite{LY1},\cite{LY2}) defined a quasi-local mass, generalizing the Brown-York 
quasi-local mass, and proved its positivity using the Jang equation. Other similar applications of the Jang equation can be found \cite{WY} and \cite{WY2} for example.\\

The problem can be stated as follows: let $(M^3,g,K)$ be an initial data set for the Einstein equation and consider the four dimensional manifold $M\times\mathbb{R}$ equipped with the {\it Riemannian} metric $\<\,,\,\>:=g\oplus dt^2$. The problem is to find a smooth function $u:M\rightarrow\mathbb{R}$ such that the hypersurface $\widehat{M}$ of $M\times\mathbb{R}$ obtained by taking the graph of $u$ over $M$, satisfies the equation
\begin{eqnarray*}
H_{\what{M}}={\rm Tr}_{\what{M}}(K)
\end{eqnarray*}
where $H_{\what{M}}$ denotes the mean curvature of $\what{M}$ in $(M\times\mathbb{R},\<\,,\,\>)$ and ${\rm Tr}_{\what{M}}(\,.\,)$ is the trace on $\what{M}$ with respect to the induced metric. This geometric problem is equivalent to solve the non-linear second order elliptic equation
\begin{eqnarray}\label{jang}
\sum_{i,j=1}^3\Big(g^{ij}-\frac{u^iu^j}{1+|\nabla u|^2}\Big)\Big(\frac{(\nabla^2u)_{ij}}{\sqrt{1+|\nabla u|^2}}-K_{ij}\Big)=0 
\end{eqnarray}
where $\nabla$ (resp. $\nabla^2$) denotes the Levi-Civita connection (resp. the Hessian) of the metric $g$, $u^i=g^{ij}u_j$ and $u_j=e_j(u)$. Note that the 
metric induced by $\<\,,\,\>$ on $\what{M}$ is
\begin{eqnarray*}
\what{g}_{ij}=g_{ij}+u_iu_j
\end{eqnarray*}
and can be viewed as a deformation of the metric $g$ on $M$. In the following, we adopt the convention that $M$ and $\what{M}$ denote respectively the Riemannian manifolds $(M,g)$ and $(M,\what{g})$. Analogously, if $\nabla$ denotes the Levi-Civita connection for $M$, then $\what{\nabla}$ denotes that on $\what{M}$ and so on. Since we assume that the initial data set $(M^3,g,K)$ comes from a spacetime satisfying the dominant energy condition, we have that the following relation holds on $\what{M}$:
\begin{eqnarray}\label{schoenyau}
 0\leq 2(\mu-|J|)\leq \what{R}-2|X|_{\what{g}}^2-2\what{\delta}(X)
\end{eqnarray}
where 
\begin{eqnarray}\label{defX}
X=\omega-\what{\nabla}\log(f),
\end{eqnarray}
$\omega$ is the tangent part of the vector field dual to $-K(\,.\,,\what{\nu})$, $f=-\<\partial_t,\what{\nu}\>$ and $\what{\nu}$ denotes the unit normal vector field to $\what{M}$ in $M\times\mathbb{R}$. All the quantities $K_{ij}$, $\mu$ and $J$ are defined on $M\times\mathbb{R}$ by parallel transport along the $\mathbb{R}$-factor. Moreover, equality occurs in (\ref{schoenyau}) if and only if $\mu=|J|$ and the second fundamental form of $\what{M}$ in 
$M\times\mathbb{R}$ is $K$. \\

It is important to note here that in  Theorem \ref{spacetime}, we assume that there is no apparent horizon in the interior of $\Omega$ so that there exists a global solution of the Jang equation which does not blow-up. 


\subsection{Proof of Theorem \ref{spacetime}}


From the work \cite{yau1} and since we assumed that $\Omega$ has no apparent horizon in its interior, there exists a smooth solution $u$ on $\Omega$  of the Jang Equation (\ref{jang}), defined  with the Dirichlet boundary condition 
\begin{eqnarray*}
 u_{|\Sigma}\equiv 0.
\end{eqnarray*}
This boundary condition ensures that the metrics $\what{g}$ and $g$ coincide on the boundary $\Sigma$ so that the Dirac operators $\D$ acting on $\sm\Sigma$ and $\what{\D}$ on $\what{\sm}\Sigma$ also coincide. Moreover, from a calculation in \cite{yau1},
we have:
\begin{eqnarray*}
\what{H}-\what{g}(X,\what{N})=f^{-1}H-\sigma|\nabla u|{\rm Tr}_\Sigma(K) 
\end{eqnarray*}
where $\what{N}$ denotes the unit outward normal vector field of $\Sigma$ in $\what{\Omega}$ and $\sigma\in\{\pm 1\}$. From this equality and since $f=-\<\partial_t,\what{\nu}\>=1/\sqrt{1+|\nabla u|^2}$, we easily see that
\begin{eqnarray}\label{boundblack}
F:=\what{H}-\what{g}(X,\what{N})\geq|\mathcal{H}|=\sqrt{H^2-{\rm Tr}_\Sigma(K)^2}.
\end{eqnarray}
Since we assume that $\Sigma$ has a spacelike mean curvature vector $\mathcal{H}$, this implies that the function $F$ is positive on $\Sigma$. From the discussion of Section \ref{jang-equation}, we also have that the resulting Riemannian manifold $\what{\Omega}$ satisfies the condition (\ref{scalarcond}) because of (\ref{schoenyau}), the vector field $X$ being defined here by (\ref{defX}). Clearly, all the assumptions of Theorem \ref{holo-princ} are fulfilled and we deduce that for all $\varphi\in\Gamma(\sm\Sigma)$:
\begin{eqnarray*}
0\le \int_\Sigma \big(\frac{1}{F}|\D \varphi|^2-\frac{F}{4}|\varphi|^2\big)\,d\Sigma,
\end{eqnarray*}
which by Inequality (\ref{boundblack}),  implies Inequality (\ref{main-holo}).

Now assume that equality is achieved. Once again we apply Theorem \ref{holo-princ} and then $\what{\Omega}$ has at least a parallel spinor field $\Phi$. In particular, $\what{\Omega}$ is Ricci flat and since it is a $3$-dimensional domain, it is flat. Moreover, if we have equality in (\ref{schoenyau}), then the second fundamental form of $\what{\Omega}$ in $M\times\mathbb{R}$ is $K_{ij}$. So we can choose a  coordinates system  $\what{x}=(\what{x}_1,\what{x}_2,\what{x}_3)$ in a neighborhood $\mathcal{U}$ of a point $p\in\Omega$ such that $\what{g}_{ij}=\delta_{ij}$. In this chart, we have:
\begin{eqnarray*}
g_{ij}=\delta_{ij}-\frac{\partial u}{\partial\what{x}_i}\frac{\partial u}{\partial\what{x}_j}
\end{eqnarray*}
and this shows that if $(\what{x}_1,\what{x_2},\what{x_3},t)$ denotes coordinates in the Minkowski spacetime, the graph of $u$ over $\mathcal{U}$ isometrically embeds in $\mathbb{R}^{3,1}$ with second fundamental form given by $K_{ij}$. Then it is clear that $\Omega$ locally embeds in the Minkowski spacetime with $K$ as second fundamental form as asserted.
\qed

As a first consequence, we have the estimate proved by the third author in \cite{Ra} for the first eigenvalue of the Dirac operator on $\Sigma$. 
\begin{corollary}\label{untrapped}
Under the same conditions of Theorem \ref{spacetime}, the first eigenvalue $\lambda_1(D_\Sigma)$ of the Dirac operator satisfies
\begin{eqnarray*}
\lambda_1(D_\Sigma)^2\geq\frac{1}{4}\inf_{\Sigma} |\mathcal{H}|^2.
\end{eqnarray*}
Moreover, if equality occurs then $\Sigma$ is connected and there exists a local isometric embedding of $\Omega$ as a spacelike hypersurface in $\mathbb{R}^{3,1}$ with $K$ as second fundamental form. 
\end{corollary}
\pf
The inequality on $\lambda_1(D_\Sigma)$ follows directly by taking $\varphi=\Phi\in\Gamma(\sm\Sigma)$ in (\ref{main-holo}) where $\Phi$ is an eigenspinor for the Dirac operator $\D$ associated with the eigenvalue $\lambda_1(\D)$ (= $\lambda_1(D_\Sigma)$). On the other hand, the second part of the equality case follows directly from Theorem \ref{spacetime}. For the connectedness of $\Sigma$, it is enough to remark that, from \cite{HMZ1}, the eigenspace associated to $\lambda_1(\D)$ corresponds to the restriction on $\Sigma$ of the space of parallel spinor fields on the domain $\what{\Omega}$ obtained by solving the Jang equation. Then, assuming that $\Sigma$ has several connected components, we fix one of them, say $\Sigma_{0}$, and define a spinor field on $\Sigma$ by
$$
\widetilde{\Phi}=\left\{
\begin{array}{ll}
\Phi_0 \qquad & \hbox{ {\rm on} } \Sigma_0 \\
0  \qquad & \hbox{ {\rm on} }\Sigma-\Sigma_0, 
\end{array}
\right.
$$
where $\Phi_0$ is an eigenspinor for the extrinsic Dirac operator $\D$ associated to the eigenvalue $\lambda_1(\D)$. It is then straightforward to check that $\widetilde{\Phi}$ is also an eigenspinor associated to $\lambda_1(\D)$ so that it comes from the restriction of a parallel spinor on $\what{\Omega}$. However since such a spinor field has constant norm, it is impossible unless $\Sigma$ is connected. 
\qed

{\textit {Proof }of Theorem \ref{liu-yau-general} :} In order to establish Inequality (\ref{liu-yau-general-Ineq}) it is sufficient to apply Inequality (\ref{main-holo}) to the restriction to $\Sigma$ of a parallel spinor field on $\mathbb R^3$. From the equality case of Theorem \ref{spacetime}, we deduce that $\Omega$ locally embeds in the Minkowski spacetime with $K$ as a second fundamental form. On the other hand, we have equality in (\ref{boundblack}) so that $\what{H}=|\mathcal{H}|$ and then equality in (\ref{liu-yau-general-Ineq}) now reads
\begin{eqnarray*}
\int_{\Sigma}\big(\what{H}-\frac{H_0^2}{\what{H}}\big)d\Sigma=0.
\end{eqnarray*}
We conclude by applying the rigidity part of Theorem $3$ in \cite{HM1} to the compact Ricci-flat manifold $\what{\Omega}$ to deduce that $\Sigma$ is connected and $|\mathcal{H}|=H_0$.
\qed


\subsection{$2$-codimensional  outer untrapped submanifolds in the Minkowski Spacetime}


In this section, we prove that Inequality (\ref{main-holo}) holds in the case of codimension two outer untrapped submanifolds of the Minkowski spacetime without any assumption on the existence of apparent horizon. More precisely, we prove:
\begin{theorem}\label{Mink}
Let $\Sigma^n$ be a codimension two outer untrapped submanifold of the $(n+2)$-dimensional Minkowski spacetime $(\mathbb{R}^{n+1,1},\<\,,\,\>)$, then Inequality $(\ref{main-holo})$ holds. Moreover, equality holds if and only if $\Sigma$ lies in a hyperplane of $\mathbb{R}^{n+1,1}$.
\end{theorem}
\pf
First we note that by assumption $\Sigma$ factorizes through a compact and connected spacelike hypersurface $\Omega$ of $\mathbb{R}^{n+1,1}$. This factorization provides  a Lorentzian orthonormal reference $\{T,N\}$ for the normal plane of $\Sigma$ in ${\mathbb{R}}^{n+1,1}$ and since $\Sigma$ is the boundary of a mean-convex domain $\Omega$ and has spacelike mean curvature vector, we deduce that the corresponding future-directed null expansions satisfy $\theta_+>0$ and $\theta_-<0$. On the other hand, from the work of Bartnik-Simon \cite{bartnik} and a straightforward generalization of Lemma $4.1$ in \cite{miaoshitam}, the submanifold $\Sigma$ spans a compact, smoothly immersed, maximal hypersurface $\Omega'$ in $\mathbb{R}^{n+1,1}$. This means that $\Sigma$ factorizes through another spacelike hypersurface $\Omega'$ of $\mathbb{R}^{n+1,1}$. The new factorization provides a different Lorentzian orthonormal reference $\{T',N'\}$ for the normal plane of $\Sigma$ in $\mathbb{R}^{n+1,1}$. In fact, it is obvious that there must be a function $f\in C^\infty(\Sigma)$, such that
$$T'=(\cosh f) T-(\sinh f) N,\qquad N'=-(\sinh f) T+(\cosh f) N.$$ 
It is clear that this new reference determines a new pair of null vectors $T'\pm N'$ and a new future-directed null expansion of ${\mathcal H}$: 
\begin{equation}\label{theta}
\theta'_+=e^f\theta_+,\qquad \theta'_-=e^{-f}\theta_-
\end{equation}
which satisfies $\theta'_+>0$ and $\theta'_-<0$. In particular, we get that $2H'=\theta'_+-\theta'_->0$. Moreover, since $\Omega'$ is maximal we have ${\rm Tr}(K')=0$ and the Gauss formula gives $R'=|K'|^2\geq 0$. Here $R'$ is the scalar curvature of $\Omega'$ equipped with the metric induced by the Minkowski spacetime and $K'$ is the associated second fundamental form. On the other hand, since $\Sigma$ has a spacelike mean curvature vector, we deduce 
\begin{eqnarray}\label{mcv}
0<|\mathcal{H}|=\sqrt{-\theta'_+\theta'_-}=\sqrt{H^{'2}-{\rm Tr}_\Sigma(K')^2}\leq H^{'}
\end{eqnarray}
so that we conclude that $\Omega'$ is such that $R'\geq 0$ and $H'>0$. Now we can apply Theorem \ref{holo-princ} to $\Omega'$ with $X\equiv 0$ and then for all $\varphi\in\Gamma(\sm\Sigma)$, we have:
\begin{eqnarray}\label{parallel-holo}
0\le \int_\Sigma \big(\frac{1}{H'}|\D \varphi|^2-\frac{1}{4}H'|\varphi|^2\big)\,d\Sigma.
\end{eqnarray}

Inequality (\ref{main-holo}) follows using  Inequality (\ref{mcv}). Assume now that equality is achieved. From the equality case of (\ref{parallel-holo}), we deduce that $\Omega'$ has at least a parallel spinor so that $\Omega'$ is Ricci flat. In particular, it has zero scalar curvature and since $R'=|K'|^2=0$, $\Omega'$ has to be totally geodesic in $\mathbb{R}^{n+1,1}$ hence $\Sigma$ lies in a hyperplane of $\mathbb{R}^{n+1,1}$. Conversely, if $\Sigma$ is $2$-codimensional  submanifold with spacelike mean curvature vector which lies in a hyperplane $\mathbb{R}^{n+1,1}$, then its second fundamental form $K$ is zero since a hyperplane $P^{n+1}$ is totally geodesic. In particular, the squared norm of the mean curvature vector of $\Sigma$ satisfies
\begin{eqnarray}\label{eq}
|\mathcal{H}|^2=H^2-{\rm Tr}_\Sigma(K) ^2=H^2,
\end{eqnarray}
where $H$ is the mean curvature of $\Sigma$ in the hyperplane $P$. Note that $|\mathcal{H}|>0$ since $H>0$. Consider now a parallel spinor field $\Phi_0$ on $\mathbb{R}^{n+1,1}$. The spinorial Gauss formula from the totally geodesic immersion of the hyperplane $P^{n+1}$ in $\mathbb{R}^{n+1,1}$ and then the one from $\Sigma^n$ into $P^{n+1}$ tell us that $\Phi_0$ satisfies for all $Y\in\Gamma(T\Sigma)$:
\begin{eqnarray*}
\nb_Y\Phi_0=-\frac{1}{2}\mult(AY)\Phi_0
\end{eqnarray*}
where $A$ is the Weingarten map of $\Sigma^n$ in $P^{n+1}$. Taking the trace of this identity gives
\begin{eqnarray*}
\D\Phi_0=\frac{1}{2}H\Phi_0=\frac{1}{2}|\mathcal{H}|\Phi_0
\end{eqnarray*}
where the last equality comes from (\ref{eq}). It is now straightforward to check that equality holds in (\ref{main-holo}) for $\varphi=\Phi_0$.
\qed

Note that Theorem \ref{liu-yau-minkowski} is obtained as a direct application of the previous result. As an application we obtain the $n$-dimensional counterpart of Corollary \ref{untrapped} in the Minkowski spacetime with an optimal rigidity statement:
\begin{corollary}\label{spacetime1}
Let $\Sigma^n$ be a codimension two outer untrapped submanifold in $\mathbb{R}^{n+1,1}$, then
\begin{equation*}\label{ineqlambda}
\left|\lambda_1(D_\Sigma)\right|\ge\frac{1}{2}\inf_\Sigma|{\mathcal H}|.
\end{equation*}
Moreover, equality occurs if and only if $\Sigma$ is a totally umbilical round sphere in a spacelike hyperplane of $\mathbb{R}^{n+1,1}$.
\end{corollary}

\pf
It is enough to apply the previous theorem to an eigenspinor for $\D$ associated with the eigenvalue $\lambda_1(\D)$ and we directly have the result. From Theorem \ref{Mink}, $\Sigma$ lies in a totally geodesic spacelike hyperplane $P^{n+1}$ with constant positive mean curvature $H$. Then the Alexandrov theorem allows to conclude that $\Sigma$ is a totally umbilical sphere in $P^{n+1}$. The converse is clear by taking the restriction of a parallel spinor of the Minkowski space to $\Sigma$ via the totally geodesic immersion of $\mathbb{R}^{n+1}$ in $\mathbb{R}^{n+1,1}$. 
\qed



\begin{thebibliography}{BHHM}

\bibitem [AR]{AR}
F. Almgren, I. Rivin, {\em The mean curvature integral is invariant under bending}, Geometry \& Topology Monographs, vol. {\bf 1}: {\em The Epstein Birthday Schrift}, 1998, 1--21.


\bibitem[AEM]{AEM}
L. Andersson, M. Eichmair and J. Metzger, \emph{Jang's equation and its
  applications to marginally trapped surfaces}, Contemporary Mathematics Complex Analysis and Dynamical Systems IV: Part 2. General Relativity, Geometry, and PDE (2011) 13--46.

\bibitem [BC]{BC} 
R. Bartnik, P. Chru{\'s}ciel, {\em Boundary value problems
for Dirac-type equations}, J. reine angew. Math., {\bf 579} (2005), 13--73. 

\bibitem [BS]{bartnik}
R. Bartnik, L. Simon, {\em Spacelike hypersurfaces with prescribed boundary values}, Comm. Math. Phys. {\bf 87} (1982),  131-152.

\bibitem [B{\"a}]{Ba2} 
C. B{\"a}r, {\em Extrinsic bounds of the Dirac operator}, 
Ann. Glob. Anal. Geom., {\bf 16} (1998), 573--596.

\bibitem [B{\"a}Ba]{BB} 
C. B{\"a}r, W. Ballmann, {\em Boundary value problems for elliptic differential operators of first order}, 
Surveys in Differential Geometry {\bf 17} (2012), 1--78.

\bibitem [BFGK]{BFGK}
H. Baum, T. Friedrich, R. Gr{\"u}newald, I. Kath,
{\it Twistor and Killing Spinors on Riemannian Manifolds},
Seminarbericht {\bf 108}, Humboldt-Universit{\"a}t zu Berlin, 1990.

\bibitem [BeS]{BeS}
R. Benedetti, R. Silhol, {\em Spin and pin$^-$ structures, immersed and embedded surfaces and a result
of Segre on real cubic surfaces}, Topology, {\bf 34} (1995), 651--678. 

\bibitem [BW]{BW} 
B. Boo\ss-Bavnbek, K.P. Wojciechowski, {\it Elliptic 
Boundary Problems for the Dirac Operator}, Birkh{\"a}user, Basel, 1993.

\bibitem [BY]{BY} 
J. Brown, J. York, {\em Quasilocal energy and conserved charges derived from the gravitational action}, Phys. Rev. D, {\bf 47} (4) (1993), 1407--1419.

\bibitem [Bu]{Bu} 
J. Bure\v s,  {\em Dirac operators on hypersurfaces}, 
Comment. Math. Univ. Carolin. {\bf 34} (1993), no. 2, 313--322. 

\bibitem [C-V]{C-V} 
S. E. Cohn-Vossen, {\em Die Verbiegung von Fl{\"a}chen im Grossen}, Fortschr. math. Wiss., {\bf 1} (1936), 33--76.

\bibitem [EMW]{EMW} 
M. Eichmair, P. Miao, X. Wang, {\em Extension of a theorem of Shi and Tam}, Calc. Var. Part. Diff. Eq., {\bf 43}, no. 1-2 (2012), 45--56.

\bibitem [HH]{HH}
J. Hass, J. Hughes, {\em Immersions of surfaces in 3-manifolds}, Topology, {\bf 24} (1985), 97--112.  

\bibitem [He]{He}
G. Herglotz, {\em {\"U}ber die Starrheit der Eifl{\"a}chen}, Abhandlungen aus dem mathematischen Seminar 
der Hansischen Universit{\"a}t, {\bf 15} (1943).

\bibitem [HM1]{HM1} 
O. Hijazi, S. Montiel, {\em A holographic principle for the existence
of parallel spinor fields and an inequality of Shi-Tam type}, 
Asian Journal of Math. {\bf 18} no.~3 (2014), 489--506.
 

\bibitem [HMZ1]{HMZ1} 
O. Hijazi, S. Montiel, X. Zhang, {\em Dirac operator on embedded hypersurfaces}, 
Math. Res. Lett, {\bf 8} (2001), 195--208.

\bibitem [HMZ2]{HMZ2} 
O. Hijazi, S. Montiel, X. Zhang, {\em Eigenvalues  of 
the Dirac operator on manifolds with boundary}, Commun. Math. Phys., 
{\bf 221} (2001), 255--265.

\bibitem [HMZ3]{HMZ3} 
O. Hijazi, S. Montiel, X. Zhang, {\em Conformal lower 
bounds for the Dirac operator of embedded hypersurfaces}, Asian J. Math., 
{\bf 6} (2002), 23--36. 

\bibitem[Ja]{jang}
P.-S. Jang, \emph{On the positivity of energy in general relativity}, J. Math.
  Phys. \textbf{19} (1978), no.~5, 1152--1155.

\bibitem [L]{L} 
G. Lam, {\em The graph cases of the Riemannian positive mass and Penrose inequalities in all dimensions}, 
thesis, Duke University, 2011.

\bibitem [LM]{LM} 
H.B. Lawson, M.L. Michelsohn, {\it Spin Geometry}, 
Princeton Math. Series, vol. 38, Princeton University Press, 1989.

\bibitem [LY1]{LY1} 
C.-C. Liu, S.-T. Yau, {\em Positivity of quasilocal mass}, 
Phys. Rev. Lett., {\bf 90} (2003), 231102--231106. 

\bibitem [LY2]{LY2} 
C.-C. Liu, S.-T. Yau, {\em Positivity of quasi-local mass II}, 
J. Amer. Math. Soc., {\bf 19} (2006), no. 1, 181--204.

\bibitem [MST]{miaoshitam}
P. Miao, Y. Shi and L.-F. Tam, \emph{On geometric problems related to
  Brown-York and Liu-Yau quasilocal mass}, Commun. Math. Phys.
  \textbf{298} (2010), no.~2, 437--459.

\bibitem [MR]{MR} 
S. Montiel, A. Ros, {\em Curves and Surfaces}, 2$^{\rm nd}$ edition, Graduate Studies in Mathematics, {\bf 69} (2009), American Mathematical Society, Rhode Island. 

\bibitem [Pi]{Pi}
U. Pinkall, {\em Regular homotopy classes of immersed surfaces}, Topology, {\bf 24} (1985), pp. 421--434.

\bibitem [R]{Ra} 
S. Raulot, {\em The Dirac operator on untrapped surfaces}, Comm. Math. Phys., {\bf 318} (2013), no.~2, 411--427. 

\bibitem [RS]{RS}
I. Rivin, J.-M. Schlenker, {\em The Schl{\"a}fli formula in Einstein manifolds with boundary}, Electronical Research Announcements 
of the AMS, {\bf 5}
(1999), 18--23. 


\bibitem[SY1]{symac}
R. Schoen and S.-T. Yau, \emph{On the proof of the positive mass conjecture in
  general relativity}, Comm. Math. Phys. \textbf{65} (1979), no.~1, 45--76.

\bibitem[SY2]{sy1}
\bysame, \emph{Proof of the positive mass Theorem II}, Comm. Math. Phys.
  \textbf{79} (1981), no.~2, 231--260.

\bibitem [ST]{ST1} Y. Shi, L.-F. Tam, {\em Positive mass theorem and the boundary
behaviors of compact manifolds with nonnegative scalar curvature}, J. Diff. Geom., 
{\bf 62} (2002), 79--125.

\bibitem [Sz]{Sz}
L.B. Szabados, {\em Quasi-local energy-momentum and angular momentum in General Relativity}, Living Rev. Relativity, {\bf 7} (2004) No 4, 1--135.

\bibitem [Tr]{Tr} A. Trautman, {\em  The Dirac operator on hypersurfaces}, 
Acta Phys. Pol., {\bf B 26} (1995), 1283--1310.

\bibitem [WY1]{WY} M.-T. Wang, S.-T. Yau, {\em A generalization of Liu-Yau's quasi-local mass}, 
Comm. Anal. Geom., {\bf 15} (2007), 249--282.

\bibitem [WY2]{WY2} M.-T. Wang, S.-T. Yau, {\em Isometric embeddings into the Minkowski space and new quasi-local mass}, Comm. Math. Phys., {\bf 288}, no.3 (2009), 919--942. 

\bibitem[We]{We} H. Weyl, {\em \"Uber die Bestimmung einer geschlossenen konvexen Fl\"ache durch ihr 
Linienelement}, Vierteljahrsschrift der naturforschenden Gesellschaft, {\bf 61} (1916), 40--72. 

\bibitem[Wi]{Wi} E. Witten, {\em A new proof of the positive energy theorem}, Commun. Math. Phys., {\bf
80} (1981), 381--402.

\bibitem[Y]{yau1}
S.-T. Yau, \emph{Geometry of three manifolds and existence of black hole due to
boundary effect}, Adv. Theor. Math. Phys. \textbf{5} (2001), no.~4, 755--767.

\end{thebibliography}
\end{document}